\documentclass[10pt]{amsart}
\usepackage{amssymb,amsmath,amsthm,amsfonts,amsopn,url}
\theoremstyle{plain}
\newtheorem{thm}{Theorem}[section]
\newtheorem*{mt*}{Main Theorem}
\newtheorem*{stfac}{Base change formulas for depth}

\newtheorem{lemma}[thm]{Lemma}
\newtheorem*{cor*}{Corollary}
\theoremstyle{definition}
\newtheorem*{definition}{Definition}
\newcommand{\ideal}[1]{\mathfrak{#1}}
\newcommand{\m}{\ideal{m}}
\newcommand{\n}{\ideal{n}}
\newcommand{\p}{\ideal{p}}
\newcommand{\q}{\ideal{q}}
\newcommand{\ffunc}[1]{\mathrm{#1}}
\newcommand{\func}[1]{\mathrm{#1} \,}
\newcommand{\depth}{\func{depth}}
\newcommand{\depthsb}{\ffunc{depth}}

\newcommand{\Ass}{\func{Ass}}
\newcommand{\pdp}[1]{\ffunc{ph. depth}_{#1}}
\newcommand{\arrow}[1]{\stackrel{#1}{\rightarrow}}

\title[Phantom depth and base change]{Phantom depth and flat base change}

\author[Neil Epstein]{Neil M. Epstein}

\address{Department of Mathematics, University of Kansas, Lawrence, Kansas  66045}

\email{epstein@math.ku.edu}

\date{Feb 10, 2005}

\thanks{The author was partially supported by the National Science Foundation.}

\keywords{tight closure, phantom depth, base change}

\subjclass[2000]{Primary 13A35; Secondary 13B40, 13C15, 13H10}

\begin{document}
\begin{abstract}
We prove that if  $f: (R,\m) \rightarrow (S,\n)$ is a flat local homomorphism, $S/\m S$ is 
Cohen-Macaulay and $F$-injective, and $R$ and $S$ share a weak test element, then a tight closure analogue
of the (standard) formula for depth and regular sequences across flat base change holds.  As a corollary,
it follows that phantom depth commutes with completion for excellent
local rings.  We give examples to show that the analogue does not hold for surjective base change.
\end{abstract}
\maketitle

All rings considered in this paper are Noetherian, local, and of positive prime characteristic $p>0$.
For such rings $R$ (among others), Hochster and Huneke \cite{HHmain} developed a theory of ``tight closure''
for finitely-generated $R$-modules.
In \cite{AbPPD}, Ian Aberbach defined a tight closure analogue of depth, called \emph{phantom depth}, and showed
that it satisfies (analogues of) many properties we expect depth to satisfy.  One such property is
a ``phantom Auslander-Buchsbaum theorem'', which is like the classical Auslander-Buchsbaum
theorem but with both depth and projective dimension\footnote{``Phantom projective dimension'',
a tight closure analogue to the classical notion of projective dimension, was introduced in \cite{HHmain}
and further developed in \cite{HHphantom}, \cite{AHH}, and \cite{nmepdep}.} replaced
by their ``phantom'' analogues.  In \cite{nmepdep} the present author showed that under mild conditions on 
$R$ and $M$, the phantom depth
of a finitely generated $R$-module $M$ is the length of \emph{any} maximal phantom
regular sequence on $M$, as Aberbach \cite{AbPPD} had proved
in the special case that $M$ has finite phantom projective dimension.

Consider the following standard, extremely useful facts:
\begin{stfac}[see e.g. \cite{BH}, section 1.2]
Let $\phi: (R,\m) \rightarrow (S,\n)$
be a flat local homomorphism of Noetherian local rings, let $M$ be a finitely generated $R$-module,
let $\mathbf{x} = x_1, \dotsc, x_a \in \m$ be an $M$-regular sequence
and let $\mathbf{y} = y_1, \dotsc, y_b \in \n$ be an $(S/\m S)$-regular sequence.
Then $\phi(\mathbf{x}), \mathbf{y}$ is an $(S \otimes_R M)$-regular sequence.  Furthermore, if $\mathbf{x}$ and $\mathbf{y}$ are \emph{maximal} regular sequences on $M$, $S/\m S$ respectively,
then the sequence $\phi(\mathbf{x}), \mathbf{y}$ is a \emph{maximal} $(S\otimes_R M)$-regular sequence.  In particular,
we have \[
\depthsb_R M + \depth S/\m S = \depthsb_S (S \otimes_R M).
\]

If instead of being flat, $\phi$ is surjective, then for any finitely-generated $S$-module $N$,
\[
\depthsb_R N = \depthsb_S N
\]
\end{stfac}
It seems natural to ask: what parts of these base change formulas hold when we replace ``depth'' and ``regular sequence''
with their phantom analogues?

Na\"ively, one would hope that all the resulting statements held verbatim.  However,
base change in tight closure theory is a messy business.  One often has to
impose conditions on the closed fiber and the residue fields, and the rings
sometimes need to share a weak test element for the proofs to work.  See \cite{HHbase},
\cite{EnflatFrat}, \cite{HaCMFI}, \cite{AbflatFreg}, \cite{AbEnTB}, \cite{HHexponent}, and \cite{SmBr}
for work along these lines, and see \cite{SiFregd} for an interesting counterexample.
These authors investigate preservation of such properties
as $F$-rationality, weak and strong $F$-regularity, whether
tight closure commutes with localization, and extension of the test ideal.

In the spirit of the some of the aforementioned papers, we prove the following ``phantom analogue''
of the flat base change formula for depth in Section~\ref{sec:fbasech}:
\begin{mt*}
Let $(R,\m) \stackrel{\phi}{\rightarrow} (S,\n)$ be a flat local homomorphism of Noetherian local rings of prime characteristic $p>0$.  Let $M$ be a finitely generated $R$-module satisfying avoidance, and suppose that $R$ and $S$ share a $q_0$-weak test element $c$ and that the closed fiber $S/\m S$ is Cohen-Macaulay and $F$-injective.  Then\footnote{Since $S/\m S$ is Cohen-Macaulay, we could have written $\dim S/\m S$,
$\pdp{S/\m S} S/\m S$, or $\pdp{S} S/\m S$ instead of $\depth S/\m S$.  However,
if the closed fiber were not Cohen-Macaulay, these four numbers could differ.
One of the problems in trying to extend the analogy to a situation where the closed fiber
is not Cohen-Macaulay would be to choose which of the above four
invariants (if any) provides the correct middle term in the displayed formula in the theorem.} \[
\pdp{R} M + \depth S / \m S = \pdp{S} (S \otimes_R M).
\]
\end{mt*}

Unfortunately, no corresponding analogue of the surjective base change formula
holds, as we show in two counterexamples in Section~\ref{sec:counter}.

\section{Background}
There are many excellent accounts of tight closure theory, including the seminal paper
\cite{HHmain} and the monograph \cite{HuTC}, so in this note we will only cover the points
most salient to our work here.  If $M$ is a finitely generated $R$-module,
where $(R,\m)$ is a Noetherian local ring of prime characteristic $p>0$, 
let $\varphi: X \rightarrow Y$ be a minimal free presentation of $M$.  If we fix
bases for the free modules $X$ and $Y$, then $\varphi$ can be thought of as
a matrix $(\varphi_{i j})$ of elements of $\m$.  
For a power $q=p^e$ of $p$, let $\varphi^q: X \rightarrow Y$ be the homomorphism 
defined by the matrix $(\varphi_{i j}^q)$.  Then we set $F^e_R(M) := F^e(M) := \func{coker} \varphi^q$.
This module is called the $q$'th \emph{Frobenius power} of $M$.
For an element $z \in M$, let $y$ be its preimage in $Y$.  Then $z^q$ is the 
image of the element $y$ in $F^e(M)$.
It is standard that $z^q$ and $F^e(M)$ are independent of the choice of free modules
and bases in the minimal free presentation, and indeed $F^e( - ) := F^e_R( - )$ can be made into
a right-exact functor from the category of finitely-generated $R$-modules
to itself.  In particular, if $f: L \rightarrow M$ is a map of finitely-generated
$R$-modules, we get a corresponding map $F^e(f): F^e(L) \rightarrow F^e(M)$.
If $i: N \hookrightarrow M$ is a submodule, then $N^{[q]}_M$ will denote
the image of $F^e(N)$ in $F^e(M)$ under the map $F^e(i)$.  Note that if $e$ and
$e'$ are positive integers and $S$ is an $R$-algebra, we have that $F^{e+e'} = F^e \circ F^{e'}$ 
and $F^e_S(S \otimes_R - ) = S \otimes_R F^e_R( - )$ as functors
on the category of finitely-generated $R$-modules.

We now have enough to define tight closure of a submodule.  Denote by 
$R^o$ the complement of the union of the minimal primes of $R$.  For a submodule
$N \subseteq M$, we say that an element $z \in M$ is in the \emph{tight closure
of $N$ in $M$} (in symbols, $z \in N^*_M$) if there is some $c \in R^o$ and some
integer $e_0$ such that for all $e \geq e_0$, $c z^{p^e} \in N^{[p^e]}_M$.
If $c$ and $q_0 = p^{e_0}$ can be chosen uniformly for all such triples $(z,M,N)$, then
we say that $c$ is a \emph{weak test element} (or a \emph{$q_0$-weak test element},
if we want to emphasize the power $q_0$) for $R$.  We say that $c$ is a
\emph{completely stable} $q_0$-weak test element for $R$ if its image is a $q_0$-weak test element
for $\hat{R}_\p$ whenever $c \in \p \in \func{Spec} \hat{R}$.  In intricate work,
Hochster and Huneke \cite{HHbase} showed that whenever $R$ is essentially of finite type
over an excellent local ring,
it has a completely stable weak test element.  In order to simplify our definitions and proofs,
we often assume that $R$ has a weak test element.

Let $G^e(M) := F^e(M) / 0^*_{F^e(M)}$, the $q$'th \emph{reduced Frobenius power} of $M$.
We say that $M$ satisfies \emph{avoidance} if for any quotient module $N$ of $M$ and any
ideal $I \subseteq R$ such that \[
I \subseteq \bigcup \bigcup_{e\geq 0} \func{Ass} G^e(N),
\] there is some $e \geq 0$ and some $\p \in \func{Ass} G^e(N)$ such that
$I \subseteq \p$.  In particular, if $\m \subseteq \bigcup \bigcup_{e\geq 0} \func{Ass}
G^e(N)$, then $\m \in \bigcup_{e \geq 0} \func{Ass} G^e(N)$.

Avoidance is a weak condition.  For example, it holds whenever
$R$ satisfies countable prime avoidance, which is the case if $R$ is complete \cite[Lemma 3]{BurchCPA}
or contains an uncountable field \cite[Remark 2.17]{HHexponent}.  It also occurs whenever the union
$\bigcup_{e \geq 0} \func{Ass} G^e(M)$ has only finitely many maximal elements, a condition for which
no counterexamples were known to exist until recently \cite{SiSwUAss}.

Next we provide the following definition of phantom $M$-regular sequences
and phantom depth.  It is \emph{a priori} different from the original one given 
in~\cite{AbPPD}, but as I show in \cite{nmepdep}, they are equivalent when $R$ has
a weak test element.

\begin{definition}
Let $R$ be a Noetherian ring of prime characteristic $p>0$ containing a weak test element, and let $M$ be a
finitely generated $R$-module.  Then we say an element $x \in R$ is
\emph{phantom $M$-regular} if $xM \neq M$ and $0:_{F^{e}(M)} x^{p^e}
\subseteq 0^{*}_{F^{e}(M)}$ for all $e \geq 0$.

A \emph{phantom zerodivisor of $M$} is an element $x \in R$ which is
not phantom $M$-regular.

A sequence $\mathbf{x} = x_{1}, \dotsc , x_{n}$ of elements of $R$ is
a \emph{phantom $M$-regular sequence} if $\mathbf{x}M \not =M$ and
$x_{i}$ is phantom $(M / (x_{1}, \dotsc , x_{i-1})M)$-regular for $1
\leq i \leq n$.

The \emph{phantom depth} of $M$ is the length of the longest phantom
$M$-regular sequence in $\m$.  It is denoted by $\pdp{\m}M$ or $\pdp{R}M$.
\end{definition}

Clearly, any $M$-regular sequence is a phantom $M$-regular sequence.  Note also that
the phantom depth of $R$ as a module over itself can be determined in a different way.
Namely, the \emph{minheight} \cite{HHphantom} of $\m$ (denoted $\func{mnht} \m$)
is defined to be $\max \{ \func{ht} \m / \p_j \mid 1 \leq j \leq t \}$,
where $\p_1, \dotsc, \p_t$ are the minimal primes of $R$.  As Aberbach notes in 
the first paragraph of the proof of \cite[Theorem 3.2.7]{AbPPD}, if $R$ is
the homomorphic image of a Cohen-Macaulay ring, then $\func{mnht} \m = \pdp{R} R$.

If $(C_., d_.)$ is a complex of finitely generated $R$-modules,
$i$ is an integer, and $Z_i = \ker d_i$ and
$B_i = \func{im} d_{i+1}$ are the cycle and boundary submodules of $C_i$, then 
we say that \emph{$H_i(C_.)$ is phantom} if $Z_i \subseteq (B_i)^*_{C_i}$ (following
\cite{HHphantom}).
The following characterization of phantom $M$-regular sequences in terms
of phantomness of Koszul homology will be crucial.
It is an analogue to the classical characterization of $M$-regular sequences
in terms of vanishing of Koszul homology:

\begin{thm}\label{thm:ghseqK} \cite{nmepdep}
Let $(R,\m)$ be a Noetherian local ring with a weak test element $c$, and
let $M$ be a finitely generated $R$-module which satisfies avoidance.  Let
$\mathbf{x} = x_1, \dotsc, x_n$ be any sequence of elements of $\m$.  Then the
following conditions are equivalent: \begin{enumerate}
\item\label{it:ghost} $\mathbf{x} = x_1, \dotsc, x_n$ is a phantom $M$-regular sequence,
\item\label{it:H1} $H_1(\mathbf{x}^{[p^e]}; F^e(M))$ is phantom for all $e \geq 0$,
\item\label{it:Hj} $H_j(\mathbf{x}^{[p^e]}; F^e(M))$ is phantom for all $e \geq 0$ and all $j \geq 1$.
\end{enumerate}
\end{thm}

In the special case where the phantom depth of a module is zero, we have the following
useful lemma:

\begin{lemma}\label{lem:phass}
Let $(R,\m)$ be a Noetherian local ring of prime characteristic $p>0$ containing a $q_0$-weak test
element $c$, and let $M$ be a 
finitely generated $R$-module.  Then the set of phantom zerodivisors for $M$ in $\m$
is the union $\bigcup \bigcup_{e\geq 0} \func{Ass} G^{e}(M)$.  Hence if $M$ satisfies
avoidance and $\pdp{R}(M) = 0$, then $\m \in \func{Ass} G^e(M)$ for some $e$.
\end{lemma}

\begin{proof}
For the first containment, suppose that $x$ be a phantom zerodivisor for $M$.  Then
there is some $e \geq 0$ such that $0 :_{F^{e}(M)} x^{q} \not
\subseteq 0^{*}_{F^{e}(M)}$.  That is, there is some $z \in F^{e}(M)
\setminus 0^{*}_{F^{e}(M)}$ with $x^{q}z = 0$.  Then $x^{q}
\overline{z} = \overline{0}$ in $G^{e}(M)$, where $\overline{z} \neq \bar{0}$,
so there is some $\p \in
\Ass G^{e}(M)$ with $x^{q} \in \p$.  Since $\p$ is prime and thus 
radical, $x \in \p$.

Conversely, let $x \in \p$ for some $\p \in \Ass G^{e}(M)$
for some $e$.  Then there is some $z \in F^{e}(M)$, $z \not \in
0^{*}_{F^{e}(M)}$, with $\p = \overline{0} :_{G^{e}(M)} \overline{z}$,
which means that $x z \in 0^{*}_{F^{e}(M)}$.  Then for all large powers
$q' \gg 0$ of $p$, \[
x^{q q'} c z^{q'} = x^{q q' - q'} \cdot (c (x z)^{q'}) = x^{q q' - q'} \cdot 0 = 0.
\]  If $x$ is phantom $M$-regular, the displayed equation shows that $c z^{q'} 
\in 0^*_{F^{e+e'}(M)}$.  Hence, $c^{q_0 + 1} z^{q'q_0} = 
c (c z^{q'})^{q_0} =0$, so that since $q'$ was any large enough
power of $p$, we conclude that $z \in 0^*_{F^e(M)}$,
contrary to assumption.  Thus, $x$ is a phantom zerodivisor for $M$.

The last statement now follows directly from the definitions.
\end{proof}

The final preliminary result that we need is the fact that the sets of associated primes
of the ``reduced Frobenius powers'' of a module are increasing:
\begin{lemma}\label{lem:Geincrease}
Let $(R,\m)$ be a Noetherian local ring of prime characteristic $p>0$ and $M$ a finitely
generated $R$-module.  Then for any $e \geq 0$, $\func{Ass} G^e(M) \subseteq
\func{Ass} G^{e+1}(M)$.
\end{lemma}

\begin{proof}
Without loss of generality, assume that $e = 0$, and let $\q \in \func{Ass} G^0(M)$.
Then there is some $z \in M \setminus 0^*_M$ such that $\q = 0^*_M : z$.  Let $I =
0^*_{F^1(M)} : z^p$.  We have \[
\q^{[p]} z^p = (\q z)^{[p]}_M \subseteq (0^*_M)^{[p]}_M \subseteq 0^*_{F^1(M)}.
\]
Hence $\q^{[p]} \subseteq 0^*_{F^1(M)} : z^p = I$.

On the other hand, let $a \in I$.  Then $a^p z^p = a^{p-1} (a z^p)
\in 0^*_{F^1(M)}$.  So for $q' \gg 0$, \[
c (a z)^{p q'} = c (a^p z^p)^{q'} = 0,
\]
which means that $az \in 0^*_M$, so $a \in 0^*_M : z = \q$.

We have shown that $\q^{[p]} \subseteq I \subseteq \q$, which means that $\q$
is minimal over $I$, so that $\q \in \func{Ass} R / I$.  Therefore there is some
$b \in R$ such that \[
\q = I : b = (0^*_{F^1(M)} : z^p) : b = 0^*_{F^1(M)} : b z^p,
\]
which proves that $\q \in \func{Ass} G^1(M)$.
\end{proof}

\section{Flat base change: proof of the main theorem}\label{sec:fbasech}
Recall that a Noetherian Cohen-Macaulay local ring $(R,\m)$ of prime characteristic $p>0$ is
said to be \emph{$F$-injective} if for any proper ideal $I$ of $R$ and any $x \in R$ such that
$x^p \in I^{[p]}$, it follows that $x \in I$.

\begin{proof}[Proof of the main theorem]
First we will prove the ``$\leq$'' direction.  Even more, it turns out that if $\mathbf{a} = a_1, \ldots, a_r$ is a phantom $M$-regular sequence in $\m$ and $\mathbf{z} = z_1, \ldots, z_s$ is an $S/\m S$ regular sequence in $\n$, then $\mathbf{\phi(a), z}$ is a phantom ${}_S (S \otimes_R M)$-regular sequence, just as one would expect from the classical case.  For brevity, let $M' = S \otimes_R M$.

First note that $\phi(\mathbf{a})$ is a phantom $M'$-regular sequence.  By induction we need only show this for the
one-element sequence $a = a_1$.  For integers $e \gg 0$, we have: \begin{align*}
0 :_{F^e_S(M')} \phi(a_1)^q  &= S \otimes_R (0 :_{F^e_R(M)} a_1^q) \\
&\subseteq S \otimes_R 0^*_{F^e_R(M)} \\
&\subseteq 0^*_{F^e_S(M')},
\end{align*}
where $q = p^e$.  The equality follows from flatness on colons.  To see the first inclusion, note first that
$0 :_{F^e_R(M)} a_1^q \subseteq 0^*_{F^e_R(M)}$ by definition of phantom $M$-regularity, and then
apply flatness of $S$ to this inclusion.  For the final inclusion, it is easy to see that the image of $S \otimes_R 0^*_{F^e_R(M)}$ 
under the map $S \otimes_R (0^*_{F^e_R(M)} \hookrightarrow F^e_R(M))$ is
contained in $0^*_{F^e_S(M')}$, and apply flatness one more time to see that the map in question 
is injective.  Since the displayed containment holds for all $e \gg 0$, it follows that $\phi(a_1)$ is a
phantom $M'$-regular element.

Since $F^e_S(S \otimes_R M) / \phi(\mathbf{a})^{[q]} F^e_S(S \otimes_R M) \cong S \otimes_R F^e_R(M / \mathbf{a}M)$ for any $e\geq 0$ and $\pdp{R}(M/\mathbf{a}M) = 0$, we can replace $M$ by $M/\mathbf{a}M$ in order to assume without loss of generality that $\pdp{R} M = 0$.  Lemma~\ref{lem:phass} then guarantees that there is some $e' \geq 0$ with $\m \in \ffunc{Ass}_R G^{e'}_R(M).$  Clearly we have that for any $e\geq 0$, $\mathbf{z}^{[q]}$ is an $(S / \m S)$-regular sequence.  Then by a standard result, e.g. \cite[Lemma 1.2.17(b)]{BH}, $\mathbf{z}^{[q]}$ is an $(S \otimes_R N)$-regular sequence for any finitely generated $R$-module $N$.  Since $F^e_S(M') = S \otimes_R F^e_R(M)$, applying this standard result to $F^e_R(M)$ together with the Koszul homology criterion for
regular sequences shows that $H_1(z_1^q, \ldots, z_s^q; F^e_S(M')) = 0$, which certainly implies phantomness, for all $e \geq 0$.  Thus, by Theorem~\ref{thm:ghseqK}, $\mathbf{z}$ is a phantom $M'$-regular sequence.

This completes one direction.  Next, assume that $\pdp{R} M = 0$, and we will prove that $\func{depth} S / \m S = \pdp{S} M'$.  Let $\mathbf{z} = z_1, \ldots, z_t$ be a maximal $(S/ \m S)$-regular sequence in $\n$.  As above, $\mathbf{z}$ is an $M'$-regular (and hence phantom $M'$-regular) sequence, so for maximality we need to show that $\pdp{S} \left(M' / \mathbf{z}M'\right) = 0$.

Since the phantom depth of $M$ is 0, we have by Lemma~\ref{lem:phass} that there is some $e \geq 0$ with $\m \in \ffunc{Ass}_R G^e_R(M)$.  This means that there is some $u \in F^e_R(M)$ with \[
\m = 0^*_{F^e_R(M)} :_R u.
\]
Note that since the sets of associated primes of the $G^{e'}_R(M)$'s are increasing (by Lemma~\ref{lem:Geincrease}), we may assume that $e \geq e_0$.  Then by flatness of $R \rightarrow S / (\mathbf{z}^{[q]})$ \cite[Lemma 1.2.17(b)]{BH} the map \[
S / (\m, \mathbf{z}^{[q]}) S = (S / \mathbf{z}^{[q]}) \otimes_R R / \m \arrow{\beta} (S / \mathbf{z}^{[q]}) \otimes_R G^e_R(M) = \frac{(S / \mathbf{z}^{[q]}) \otimes F^e_R(M)}{(S / \mathbf{z}^{[q]}) \otimes 0^*_{F^e_R(M)}}
\]
which sends 1 to $u$ is an injection

Also, setting $q = p^e$, there is some $b \in S$ such that \[
\n = (\m, \mathbf{z}^{[q]})S :_S b,
\]
since $\mathbf{z}^{[q]}$ is a maximal $(S / \m S)$-regular sequence in $\n$.
We have then that the map \[
S / \n \arrow{\alpha} S / (\m, \mathbf{z}^{[q]}) S
\]
that sends 1 to $b$ is an injection as well.

Now, consider the following
\begin{lemma}\label{lem:CMFIflat}
Let $(R,\m) \rightarrow (S,\n)$ be a flat local homomorphism of Noetherian local rings, both of prime characteristic $p>0$, with Cohen-Macaulay $F$-injective closed fiber and a shared 
$q_0$-weak test element $c$.  Suppose $\mathbf{z}$ is a system of parameters for the $S$-module $(S /\m S)$ and the image of $b$ is nonzero in $S 
/ (\m, \mathbf{z}) S.$  If $N$ is a finitely generated $R$-module and $u$ is not in $0^*_N$, then $bu$ is not in $0^*_{(S / \mathbf{z}) \otimes_R N}$, where the tight closure is taken over $S$.
\end{lemma}
This lemma has the same conclusion as \cite[Lemma 3.1]{AbEnTB}, and it has exactly the same proof (except that ``for all $q$'' must be replaced by ``for all $q \gg 0$''), although the hypotheses differ.  For completeness, we reproduce a version of the proof here:

\begin{proof}[Proof of Lemma~\ref{lem:CMFIflat}]
We prove the contrapositive.  That is, assuming that $bu \in 0^*_{(S / \mathbf{z}) \otimes_R N}$ (with the tight closure
taken over $S$), we will
show that $u \in 0^*_N$.

We have that for all powers $q \geq q_0$ of $p$, \[
\overline{c(bu)^q} = \bar{0} \in S / \mathbf{z}^{[q]} \otimes_R F^e_R(N).
\]  Then by flatness of the map $R \rightarrow S / \mathbf{z}^{[q]}$, \begin{equation}\label{eq:*}
\overline{b^q} \in 0 :_{S / \mathbf{z}^{[q]}} \overline{c u^q} = S / \mathbf{z}^{[q]} \otimes_R (0 :_R c u^q).
\end{equation}
If $c u^q \neq 0$, then $0 :_R c u^q \subseteq \m$, from which we conclude that \begin{equation}\label{eq:**}
S / \mathbf{z}^{[q]} \otimes_R (0 :_R c u^q) \subseteq S / \mathbf{z}^{[q]} \otimes_R \m = \frac{(\m, \mathbf{z}^{[q]}) S}{(\mathbf{z}^{[q]})}.
\end{equation}
Combining (\ref{eq:*}) with (\ref{eq:**}), we have that $b^q \in (\m, \mathbf{z}^{[q]})S$, from which we conclude, by $F$-injectivity of $S / \m S$, that $b \in (\m, \mathbf{z})S$, which contradicts our assumption on $b$.  Hence $c u^q = 0$.  Since $q$ was allowed to be any power of $p$ larger than $q_0$, it follows that $u \in 0^*_N$.
\end{proof}

It follows from Lemma~\ref{lem:CMFIflat} that, in our case, \[
bu \notin 0^*_{(S / \mathbf{z}^{[q]}) \otimes_R F^e_R(M)} = 0^*_{F^e_S((S / \mathbf{z}) \otimes_R M)}
= 0^*_{F^e_S(M'/ \mathbf{z}M')},
\]
where the tight closures are computed over $S$.  Thus, $S/ \n$ injects into $G^e_S\left(M' / \mathbf{z}M'\right)$, so that $\pdp{S} (M'/ \mathbf{z}M')= 0$, which means that $\mathbf{z}$ is indeed maximal as a phantom $M'$-regular sequence.

Finally, consider the case where $\pdp{R}M >0$.  Then if $\mathbf{a}$ is a maximal phantom $M$-regular sequence and $\mathbf{z}$ is a maximal $(S/\m S)$-regular sequence, then $M / \mathbf{a}M$ has phantom depth 0, so we can apply the above to show that $\phi(\mathbf{a}), \mathbf{z}$ is a maximal phantom $M'$-regular sequence.
\end{proof}

\begin{cor*}[Phantom depth is unaffected by completion]
Suppose $R$ is a Noetherian local ring of prime characteristic $p>0$ containing a completely stable weak test element.  For any finitely generated $R$-module $M$ that satisfies avoidance, \[
\pdp{R}M = \pdp{\hat{R}} \hat{M}.
\]
\end{cor*}

\section{Surjective base change: counterexamples}\label{sec:counter}
In general, if $(R,\m)$ is a Noetherian local ring, $S = R/I$ is a quotient of it, and $M$
is a finitely generated $S$-module, we have that $\ffunc{depth}_{\m} M =
\ffunc{depth}_{\m/I} M$.  However, unlike depth, phantom depth can depend on the ring over which it is 
calculated.  In particular, $\pdp{R} M$ may differ from $\pdp{R/I} M$.

For instance, consider the following situation from \cite[Remark 2.7]{HHphantom}, which was
also considered by Aberbach in \cite{AbPPD}.  Let 
$k$ be a field of prime characteristic $p>0$, let $T = k[[Y_1, \dotsc, Y_n, Z]]$ where $n > 1$,
let $J$ be the ideal $(Y_1 Z, \dotsc, Y_n Z)$ of $T$, set $R = T / J$ and $\m = \m_R$,
and let the images of each $Y_j$ be denoted by $y_j$, and the image of $Z$ by $z$.
Then put $x = z - y_1$ and $I = (x)$.  We have that
$x$ is a 
nonzerodivisor of $R$, so it is certainly a phantom
$R$-regular element.  Hence, $\ffunc{ppd}_R(R / xR) = 1$, so by Aberbach's phantom
Auslander-Buchsbaum theorem \cite[Theorem 3.2.7]{AbPPD}, \[
\pdp{R} R / I = \func{mnht} \m - \ffunc{ppd}_R (R / xR) = 1 -1 = 0.
\]
On the other hand, \[
\pdp{R/I} R/I = \func{mnht} \m / I = n-1 > 0 = \pdp{R} R/I.
\]

In some sense, surjective base change for phantom depth fails in the above example because $R$ is not equidimensional.  However,
surjective base change may fail even with equidimensional rings.  Let $(R,\m)$ be a reduced Noetherian 
local ring of 
characteristic $p>0$ which contains a test element $c$ and 
an ideal $I$ such that $\m$ is not minimal over $I$, but such that $\m$ 
is an associated 
prime of $I^*$, and assume further that the ring $R/I$ is equidimensional.  

Then \[
\pdp{R/I} R/I = \func{mnht} \m/I = \func{ht} \m/I > 0.
\]
On the other hand, since
$\m$ is associated to $I^*$, there is some $z \in R \setminus I^*$ such that $\m = I^* :_R z$.  Now 
suppose that there is some $a \in \m$ which is phantom ${}_R(R/I)$-regular.  Since
$a \in \m = I^* : z$, $az \in I^*$, so that for any power $q$ of $p$, \[
c a^q z^q \in I^{[q]}.
\]
Since $a$ is phantom ${}_R(R/I)$-regular, this implies that $c z^q \in 
\left(I^{[q]}\right)^*$, so that $c^2 z^q \in I^{[q]}$.  Since this holds for 
all $q$, it follows that $z \in I^*$, which is a contradiction.

Thus, $\m$ has no phantom ${}_R(R/I)$-regular elements.  That is, \[
\pdp{R} R/I = 0.
\]

For concreteness, note that such a ring $R$ and ideal $I$ is given (and proved to be such, 
apart from the equidimensionality and test element hypotheses) in 
\cite[Example 2.13]{HHexponent}.  They let \[
R = K[X,Y,U,V] / (X^3 Y^3 + U^3 + V^3) = K[x,y,u,v],
\]
where $K$ is a field of prime characteristic $p$, where $p \neq 0, 3$.  They
let \[
I = (u, v, x^3) \subseteq R.
\]
The only parts of our criteria not explicitly stated by Hochster and Huneke in \cite{HHexponent} are the
existence of a test element $c$ and the equidimensionality of $R/I$.  But as they do state, 
$R/I \cong K[x,y] / (x^3)$, which is clearly Cohen-Macaulay, hence equidimensional.  Note also that $R$ itself
is equidimensional, since it is an integral domain.  Furthermore, $R$ is a 
finitely generated algebra over a field, which implies that it has a completely
stable test element.

One might protest at this point that $R$ is not a local ring.  However, we may replace $R$ by
$R_\m$ and $I$ by $I R_\m$, where $\m = (x,y,u,v)$
is the homogeneous graded maximal ideal, without affecting any of the criteria we
needed.

\section*{Acknowledgements}
I wish to thank my advisor, Craig Huneke, for many helpful conversations and constant
encouragement.  Thanks also to Ian Aberbach for teaching me how to use fibers.

\providecommand{\bysame}{\leavevmode\hbox to3em{\hrulefill}\thinspace}
\providecommand{\MR}{\relax\ifhmode\unskip\space\fi MR }
\providecommand{\MRhref}[2]{%
  \href{http://www.ams.org/mathscinet-getitem?mr=#1}{#2}
}
\providecommand{\href}[2]{#2}

\end{document}